\documentclass[reqno]{amsart}
\usepackage{amsthm,amsmath,amssymb}
\usepackage{stmaryrd,mathrsfs}
\usepackage{esint}

\newcommand{\ud}[0]{\,\mathrm{d}}

\newcommand{\ceil}[1]{\lceil #1 \rceil}

\newcommand{\dist}[0]{\operatorname{dist}}
\newcommand{\sinc}[0]{\operatorname{sinc}}

\newcommand{\abs}[1]{|#1|}

\newcommand{\Norm}[2]{\|#1\|_{#2}}

\newcommand{\BNorm}[2]{\Big\|#1\Big\|_{#2}}
\newcommand{\pair}[2]{\langle #1,#2 \rangle}

\newcommand{\ave}[1]{\langle #1\rangle}

\newcommand{\bddlin}[0]{\mathscr{L}}

\newcommand{\supp}[0]{\operatorname{supp}}
\newcommand{\R}{\mathbb{R}}

\newcommand{\Z}{\mathbb{Z}}
\newcommand{\T}{\mathbb{T}}

\newcommand{\prob}[0]{\mathbb{P}}
\newcommand{\Exp}[0]{\mathbb{E}}
\newcommand{\D}[0]{\mathbb{D}}

\newcommand{\radem}[0]{\varepsilon}

\newcommand{\good}[0]{\operatorname{good}}
\newcommand{\bad}[0]{\operatorname{bad}}

\numberwithin{equation}{section}

\makeatletter
  \let\c@equation\c@subsection
\makeatother

\theoremstyle{plain}

\newtheorem{theorem}[subsection]{Theorem}

\newtheorem{corollary}[subsection]{Corollary}

\theoremstyle{definition}

\theoremstyle{remark}

\begin{document}

\title[Vector-valued singular integrals]%[Analysis in UMD spaces]%{A new proof of the vector-valued $T(1)$ theorem}
%{A simplified approach to two cornerstones of analysis in UMD spaces}
{Foundations of vector-valued singular integrals revisited---with random dyadic cubes}

\author[T.~P.\ Hyt\"onen]{Tuomas P.\ Hyt\"onen}
\address{Department of Mathematics and Statistics, University of Helsinki, Gustaf H\"all\-str\"omin katu 2b, FI-00014 Helsinki, Finland}
\email{tuomas.hytonen@helsinki.fi}

\date{\today}

% 42B20 Singular and oscillatory integrals (Calder\'on-Zygmund, etc.) 
% 46B09 Probabilistic methods in Banach space theory
% 46E40 Spaces of vector- and operator-valued functions
% 47A60 Functional calculus
% 47F05 Partial differential operators
% 60G46 Martingales and classical analysis

\maketitle

\begin{abstract}
The vector-valued $T(1)$ theorem due to Figiel, and a certain square function estimate of Bourgain for translations of functions with a limited frequency spectrum, are two cornerstones of harmonic analysis in UMD spaces. In this paper, a simplified approach to these results is presented, exploiting Nazarov, Treil and Volberg's method of random dyadic cubes, which allows to circumvent the most subtle parts of the original arguments.
\end{abstract}

\section{Introduction}

The investigation, during the 1980's \cite{Bourgain:86,Burkholder:83,Figiel:89,McConnell:84}, of the interrelation between the boundedness properties of vector-valued singular integral operators, and geometric or probabilistic properties of the underlying Banach space, culminated in the end of the decade in the proof of the full $T(1)$ theorem for UMD space -valued functions by Figiel \cite{Figiel:90}. This result remains a heavy piece of hard analysis, whose proof depends on subtle combinatorial arguments related to the rearrangements of dyadic cubes~\cite{Figiel:89}. An alternative Fourier-analytic proof of Weis and the author \cite{HytWei:06}, in turn, relies on a delicate square function estimate of Bourgain \cite{Bourgain:86} concerning the action of translation operators on functions of bounded frequency spectrum, whose proof, which predates and foreshadows the  rearrangements of Figiel, is similar in spirit and on the same level of difficulty.

Interestingly, some of the complications in Figiel's proof are related to essentially similar configurations of cubes (a smaller cube close to the boundary of a much larger one), which were termed ``bad'' by Nazarov, Treil and Volberg in their $T(1)$ theorem for non-homogeneous spaces \cite{NTV:97,NTV:03}, and which they overcame by a probabilistic argument using a random choice of a dyadic system instead of a fixed one. This suggests the possibility of simplifying Figiel's proof with the help of the random dyadic systems, and indeed this idea works out surprisingly nicely in this paper. 

% The notion of good and bad cubes to be used here is slightly different from the original one, involving only one random dyadic system rather than two independent copies, and may turn out to be useful in other contexts. In fact, this idea grew out of our investigation of some weighted inequalities by the related techniques (work in progress).

Since there exists already a reasonably streamlined vector-valued argument for the paraproduct operators (see Figiel and Wojtaszczyk \cite{FigWoj:01}, who attribute it to Bourgain), the present article concentrates on a short proof of the following special $T(1)$ theorem, whose quantitative form seems also new as such:

\begin{theorem}[Figiel \cite{Figiel:90}]\label{thm:Figiel}
Let $X$ be a UMD space, $p\in(1,\infty)$, and $\beta_{p,X}$ be the unconditionality constant of martingale differences in $L^p(\R^n;X)$. 
Let $T$ be a Calder\'on--Zygmund operator on $\R^n$ which satisfies the standard kernel estimates, the weak boundedness property $\abs{\pair{1_I}{T1_I}}\leq C\abs{I}$ for all cubes $I$, and the vanishing paraproduct conditions $T(1)=T^*(1)=0$.
Then $T$ extends to a bounded linear operator on $L^p(\R^n;X)$, and more precisely
\begin{equation*}
  \Norm{T}{\bddlin(L^p(\R^n;X))}
  \leq C\beta_{p,X}^2,
\end{equation*}
where $C$ depends only on the dimension $n$ and the constants in the standard estimates and the weak boundedness property.
\end{theorem}

It would be natural to conjecture that the correct estimate should be $C\beta_{p,X}$, but the quadratic bound in terms of $\beta_{p,X}$ is the best that is known in general UMD spaces even for the Hilbert transform. In the scalar-valued case, both singular integral and martingale transform norms in $L^p$ grow like $\max\{p,(p-1)^{-1}\}$, and a similar behaviour in noncommutative $L^p$ spaces has been verified by Randrianantoanina \cite{Randria:02} for martingale transforms, and recently by Parcet \cite{Parcet:09} for singular integrals, so the conjecture is true in these special cases.

In abstract UMD spaces,  a linear bound has only been shown for some special classes of operators with an even kernel (see \cite{GMS:10,Hyt:07}). The largest class of operators for which the quadratic bound was previously known seems to consist of the Fourier multiplier operators with symbol $\sigma$ satisfying $\abs{\partial^{\alpha}\sigma(\xi)}\leq C\abs{\xi}^{-\abs{\alpha}}$ for all multi-indices $\abs{\alpha}\leq n+1$. This estimate can be extracted out of the proof of McConnell~\cite{McConnell:84}, although it is not explicitly formulated in his paper. The above multiplier condition, which is stronger than the usual Mihlin or H\"ormander conditions, implies the standard estimates for the corresponding convolution kernel, so these operators also fall under the scope of the above theorem. For more general multiplier classes with fewer derivatives, the known estimates give higher powers of~$\beta_{p,X}$.

As another application of random dyadic systems to vector-valued analysis, I also provide a simpler proof of the mentioned square function estimate of Bourgain (see Theorem~\ref{thm:Bourgain}). Besides the Fourier-analytic approach to the vector-valued $T(1)$ theorem, Bourgain's inequality is also central to various other results in harmonic analysis in UMD spaces, so it seems useful to make it more approachable.

It should be emphasized that this paper is largely expository in character, and even the simplified proofs borrow their general structure and much of the details from the original arguments of Bourgain \cite{Bourgain:86} and Figiel \cite{Figiel:89,Figiel:90}. For Theorem~\ref{thm:Figiel}, a completely selfcontained proof will be provided in order to convince the reader that it can indeed be done in just about five pages. (The proof has not yet started, and it will be finished on page~\pageref{p:endFigiel}.) In the case of Bourgain's inequality, a couple of simpler lemmas (which are important on their own, and quite well known among experts in vector-valued harmonic analysis) from his original paper \cite{Bourgain:86} are taken for granted here. The principal novelty in both proofs consists of avoiding the more subtle points with the help of the probabilistic argument of Nazarov, Treil and Volberg; everything else is basically due to the original authors, even when this is not indicated explicitly at every step.

\subsection*{Acknowledgement}
Support by the Academy of Finland, grants 130166, 133264 and 218148, is gratefully acknowledged.
Some of the ideas in this paper grew out of my discussions with Michael Lacey.

\section{Random dyadic systems; good and bad cubes}

The following construction of random dyadic systems is, up to some details, from Nazarov, Treil and Volberg \cite{NTV:97,NTV:03}, these details being as in \cite{Hyt:A2}. Only one random system, rather than two, will be used here; this is akin to \cite{Hyt:A2,NTV:97} but in contrast to the probably best-known appearance of this kind of constructions in \cite{NTV:03}.

Let $\mathscr{D}^0:=\bigcup_{j\in\Z}\mathscr{D}^0_j$, $\mathscr{D}^0_j:=\{2^{-j}([0,1)^n+m):m\in\Z^n\}$ be the standard system of dyadic cubes. For every $\beta=(\beta_j)_{j\in\Z}\in(\{0,1\}^n)^{\Z}$, consider the dyadic system $\mathscr{D}^{\beta}:=\bigcup_{j\in\Z}\mathscr{D}^{\beta}_j$, where $\mathscr{D}^{\beta}_j:=\mathscr{D}^0_j+\sum_{i>j}2^{-i}\beta_i$. It is also convenient to define the shift of an individual cube $I\in\mathscr{D}^0$ by the formal shift parameter $\beta$ by using the same truncation procedure, $I+\beta:=I+\sum_{i:2^{-i}<\ell(I)}2^{-i}\beta_i$, so that $\mathscr{D}^{\beta}=\{I+\beta:I\in\mathscr{D}^0\}$.

On the space $(\{0,1\}^n)^{\Z}$, consider the natural probability $\prob_{\beta}$, which makes the coordinates $\beta_j$ independent and uniformly distributed over the set $\{0,1\}^n$. This induces a probability on the family of all dyadic systems $\mathscr{D}^{\beta}$ as defined above.

Consider for the moment a fixed dyadic system $\mathscr{D}=\mathscr{D}^{\beta}$ for some $\beta$. A cube $I\in\mathscr{D}$ is called \emph{bad} (with parameters $r\in\Z_+$ and $\gamma\in(0,1)$) if there holds
\begin{equation*}
  \dist(I,J^c)\leq \ell(I)^{\gamma}\ell(J)^{1-\gamma}\qquad
  \text{for some }J=I^{(k)},\quad k\geq r,
\end{equation*}
where $\ell(I)$ denotes the sidelength and $I^{(k)}$ the $k$th dyadic ancestor of $I$. Otherwise, $I$ is said to be \emph{good}.

Fixing a $I\in\mathscr{D}^0$, consider the random event that its shift $I+\beta$ is bad in $\mathscr{D}^{\beta}$. The badness only depends on the relative position of $I+\beta$ inside the bigger cubes in $\mathscr{D}^{\beta}$, which is determined by the coordinates $\beta_j$ with $2^{-j}\geq\ell(I)$. On the other hand, the absolute position of $I+\beta$ only depends on the coordinates $\beta_j$ with $2^{-j}<\ell(I)$, and hence the badness and position of $I+\beta$ are independent under the random choice of $\beta$. For reasons of symmetry it is obvious that the probability $\prob_{\beta}(I+\beta\text{ is bad})$ is independent of the cube $I$, and we denote it by $\pi_{\bad}$; similarly one defines $\pi_{\good}=1-\pi_{\bad}$. It is easy to see that the probability of $(I+\beta)^{(k)}$ making $I+\beta$ bad for a fixed $k$ (an event which depends on $\sum_{j:\ell(I)\leq 2^{-j}<2^k\ell(I)}\beta_j 2^{-j}$) is
\begin{equation*}
  \big(1-2\cdot 2^{-k}\ceil{2^{k(1-\gamma)}}\big)^n\leq
  2n\cdot 2^{-k}(2^{k(1-\gamma)}+1)\leq 4n 2^{-k\gamma},
\end{equation*}
and hence
\begin{equation*}
  \pi_{\bad}
  \leq\sum_{k=r}^{\infty}4n 2^{-k\gamma}
  =\frac{4n 2^{-r\gamma}}{1-2^{-\gamma}}.
\end{equation*}
The only thing that is needed about this number in the present paper, as in \cite{Hyt:A2}, is that $\pi_{\bad}<1$, and hence $\pi_{\good}>0$, as soon as $r$ is chosen sufficiently large. We henceforth consider the parameters $\gamma$ and $r$ being fixed in such a way.

\section{The dyadic representation of an operator}

For $\mathscr{D}=\mathscr{D}^{\beta}$, let $\Exp_j$ denote the conditional expectation with respect to $\mathscr{D}_j$, and $\D_j:=\Exp_{j+1}-\Exp_j$. These operators can be conveniently represented in terms of the Haar functions $h_I^{\theta}$, $\theta\in\{0,1\}^n$, defined as follows: For $n=1$,
\begin{equation*}
  h^0_I:=\abs{I}^{-1/2}1_I,\qquad
  h^1_I:=\abs{I}^{-1/2}(1_{I_+}-1_{I_-}),
\end{equation*}
where $I_+$ and $I_-$ are the left and right halves of $I$, and in general,
\begin{equation*}
  h^{\theta}_I(x)=h^{(\theta_1,\ldots,\theta_n)}_{I_1\times\cdots\times I_n}(x_1,\ldots,x_n)
  :=\prod_{i=1}^n h^{\theta_i}_{I_i}(x_i).
\end{equation*}
Then
\begin{equation*}
  \Exp_j f=\sum_{I\in\mathscr{D}_j}h^0_I\pair{h^0_I}{f},\qquad
  \D_j f=\sum_{I\in\mathscr{D}_j}\sum_{\theta\in\{0,1\}^n\setminus\{0\}}
      h^{\theta}_I\pair{h^{\theta}_I}{f}.
\end{equation*}
A frequently occurring object is the translation of a dyadic cube $I$ by $m\in\Z^n$ times its sidelength $\ell(I)$; this will be abbreviated as $I\dot+m:=I+\ell(I)m$.

The convergence of $\Exp_j f$ to $f$ as $j\to\infty$ and to $0$ as $j\to-\infty$ (both a.e. and in $L^p(\R^n)$ for $p\in(1,\infty)$) leads to Figiel's representation of an operator $T$ as the telescopic series
\begin{align*}
  \pair{g}{Tf}
  &=\sum_{j\in\Z}\big(\pair{\Exp_{j+1} g}{T\Exp_{j+1} f}
       -\pair{\Exp_j g}{T\Exp_j f}\big) \\
  &=\sum_{j\in\Z}\big(
    \pair{\D_j g}{T\D_j f}
    +\pair{\Exp_j g}{T\D_j f}+\pair{\D_j g}{T\Exp_j f}\big)
    =: A+B+C,
\end{align*}
where, upon expanding in terms of the Haar functions,
\begin{align*}
  A &=\sum_{m\in\Z^n}\sum_{I\in\mathscr{D}}
    \sum_{\eta,\theta\in\{0,1\}^n\setminus\{0\}}
     \pair{g}{h_{I\dot+m}^{\eta}}\pair{h_{I\dot+m}^{\eta}}{Th_I^{\theta}}
      \pair{h_I^{\theta}}{f}, \\
  B &=\sum_{m\in\Z^n}\sum_{I\in\mathscr{D}}
    \sum_{\theta\in\{0,1\}^n\setminus\{0\}}
     \pair{g}{h_{I\dot+m}^{0}}\pair{h_{I\dot+m}^{0}}{Th_I^{\theta}}
      \pair{h_I^{\theta}}{f} \\
  &=\sum_{m\in\Z^n}\sum_{I\in\mathscr{D}}
  \sum_{\theta\in\{0,1\}^n\setminus\{0\}}
     \pair{g}{h_{I\dot+m}^{0}-h_I^0}\pair{h_{I\dot+m}^{0}}{Th_I^{\theta}}
      \pair{h_I^{\theta}}{f} \\
  &\qquad+\sum_{I\in\mathscr{D}}\sum_{\theta\in\{0,1\}^n\setminus\{0\}}
     \ave{g}_I\pair{T^*1}{h_I^{\theta}}\pair{h_I^{\theta}}{f}=:B^0+P,
\end{align*}
and the term $C$ is essentially dual to $B$ and can be treated similarly by splitting into a cancellative part $C^0$ and a paraproduct part $Q$. It is quite explicit in the above formula that $P$ vanishes under the condition that $T^*1=0$, and similarly $Q$ is zero if $T1=0$.

If $T$ is a Calder\'on--Zygmund singular integral
\begin{equation*}
  Tf(x)=\int_{\R^n}K(x,y)f(y)\ud y,\qquad x\notin\supp f,
\end{equation*}
which satisfies the standard estimates $\abs{K(x,y)}\leq C\abs{x-y}^{-n}$ and
\begin{equation*}
  \abs{K(x+h,y)-K(x,y)}+\abs{K(x,y+h)-K(x,y)}
  \leq\frac{C\abs{h}^{\alpha}}{\abs{x-y}^{n+\alpha}}
\end{equation*}
for $\abs{x-y}>2\abs{h}$, as well as the weak boundedness property $\abs{\pair{1_I}{T1_I}}\leq C\abs{I}$ for all cubes $I$, then the Haar coefficients of $T$ satisfy
\begin{align*}
  \abs{\pair{h_{I\dot+m}^{\eta}}{Th_I^{\theta}}}
  \lesssim (1+\abs{m})^{-n-\alpha},
  \qquad(\eta,\theta)\in\{0,1\}^{2n}\setminus\{(0,0)\}.
\end{align*}
(Here and below, the notation $\lesssim$ indicates an inequality of the type ``$\leq C\times\ldots$'', where $C$ may depend at most on the dimension $n$ and the constants appearing in the Calder\'on--Zygmund conditions.) The above estimate was observed by Figiel, and it follows by elementary computations, using the H\"older estimate for the kernel when $m\notin\{-1,0,1\}^n$, the pointwise bound when $m\in\{-1,0,1\}^n\setminus\{0\}$, and the pointwise bound in combination with the weak boundedness property for $m=0$. Then it is easy to check that all the above expansions converge absolutely for example for $f\in C_c^1(\R^n;X)$ and $g\in C_c^1(\R^n;X^*)$.

Now the above expansions of $\pair{g}{Tf}$ are valid with any dyadic system $\mathscr{D}=\mathscr{D}^{\beta}$. Hence they are also valid if we take the average over a random choice of the dyadic system. For the manipulation of such averages, it is convenient to organise the summations over the fixed reference system of dyadic cubes $\mathscr{D}^0$, so that the summation condition does not depend on any random variable.

A basic observation is the following. Let
\begin{align*}
  \pi_{\good}:=\prob_{\beta}(I+\beta\text{ is good})=\Exp_{\beta}1_{\good}(I+\beta),
\end{align*}
recalling that this number is independent of the particular cube $I$. For any function $\phi(I)$ defined on the collection of all cubes $I$, we then have
\begin{align*}
   \pi_{\good}\Exp_{\beta}\sum_{I\in\mathscr{D}^{\beta}}\phi(I)
   &=\sum_{I\in\mathscr{D}^{0}}
      \Exp_{\beta}1_{\good}(I+\beta)\Exp_{\beta}\phi(I+\beta) \\
    &=\sum_{I\in\mathscr{D}^{0}}
      \Exp_{\beta}\big(1_{\good}(I+\beta)\phi(I+\beta)\big) 
    =\Exp_{\beta}\sum_{I\in\mathscr{D}^{\beta}_{\good}}\phi(I+\beta).
\end{align*}
The second equality used the crucial fact that the event that $I+\beta$ is good is independent of the position of the cube $I+\beta$, and hence of the function $\phi(I+\beta)$. The conclusion of this computation is this: inside the average over the random choice of our dyadic system, any summation over the dyadic cubes may be restricted to the good ones, and the final result is only changed by the absolute multiplicative factor $\pi_{\good}$. This gives rise to the final form of our dyadic representation,
\begin{align*}
  \pair{g}{Tf}
  =\frac{1}{\pi_{\good}}\Exp_{\beta}(A_{\good}+B^0_{\good}+C^0_{\good})
    +\Exp_{\beta}(P+Q),
\end{align*}
where e.g.
\begin{align*}
  A_{\good}=A_{\good}^{\beta}
  =\sum_{m\in\Z^n}\sum_{I\in\mathscr{D}_{\good}^{\beta}}
    \sum_{\eta,\theta\in\{0,1\}^n\setminus\{0\}}
     \pair{g}{h_{I\dot+m}^{\eta}}\pair{h_{I\dot+m}^{\eta}}{Th_I^{\theta}}
      \pair{h_I^{\theta}}{f},
\end{align*}
and the terms $B^0_{\good}$ and $C^0_{\good}$ are defined in an analogous manner by restricting the summations over $I\in\mathscr{D}$ appearing in $B^0$ and $C^0$ to $I\in\mathscr{D}_{\good}$ only. This could have been done for the paraproduct terms as well, but the known arguments for handling them do not seem to gain any particular simplification from this reduction.

\section{Estimating the expansions as martingale transforms}

The estimation of the operator norm of $T$ via the size of the pairings $\pair{g}{Tf}$ has now been reduced to the estimation of the cancellative parts $A_{\good}$, $B^0_{\good}$ and $C^0_{\good}$ as well as, in general, the paraproducts $P$ and $Q$ which are assumed to be zero here. The randomisation over the choice of the dyadic system was already fully exploited in making this reduction, and the remaining estimates will be carried out uniformly for any fixed choice of $\mathscr{D}=\mathscr{D}^{\beta}$.

Figiel's key idea for the estimation of $A$ is the interpretation of the shifted functions $h^{\eta}_{I\dot+m}$ as martingale transforms of the $h^{\theta}_I$, for each $m\in\Z^n$.
For a single $I\in\mathscr{D}$, this is easily achieved by defining the two-element martingale difference sequence
\begin{align*}
  d_{I,m,u}^{\eta,\theta}
  :=\frac{1}{2}\big(h_I^{\eta}+(-1)^u h_{I\dot+m}^{\theta}\big),\qquad u=0,1,
\end{align*}
so that
\begin{align*}
  h_I^{\eta}=d_{I,m,0}^{\eta,\theta}+d_{I,m,1}^{\eta,\theta},\qquad
  h_{I\dot+m}^{\theta}=d_{I,m,0}^{\eta,\theta}-d_{I,m,1}^{\eta,\theta}.
\end{align*}
A little more tricky is to do this in such a way that the $d_{I,m,u}^{\eta,\theta}$ still form a martingale difference sequence when also the cube $I$ is allowed to vary. This is not true for all $I\in\mathscr{D}_{\good}$, but can be achieved for appropriate subcollections which partition $\mathscr{D}_{\good}$.

For each $m$, let $M=M(m):=\max\{r,\ceil{(1-\gamma)^{-1}\log_2^+\abs{m}}\}$. Let then $a(I):=\log_2\ell(I)\mod M+1$, and define $b(I)$ to be alternatingly $0$ and $1$ along each orbit of the permutation $I\mapsto I\dot+m$ of $\mathscr{D}$. We claim that if $(a(I),b(I))=(a(J),b(J))$ for two different cubes $I,J\in\mathscr{D}_{\good}$, then the cubes satisfy the following \emph{$m$-compatibility} condition: either the sets $I\cup(I\dot+m)$ and $J\cup(J\dot+m)$ are disjoint, or one of them, say $I\cup(I\dot+m)$, is contained in a dyadic subcube of $J$ or $J\dot+m$.

Suppose first that $\ell(I)=\ell(J)$. Then $b(I)=b(J)$ ensures that $I\dot+m\neq J$ and $J\dot+m\neq I$, so the disjointness condition holds. Let then for example $\ell(I)<\ell(J)$. Then $a(I)=a(J)$ implies that in fact $\ell(J)\geq 2^{M+1}\ell(I)$. If $I$ intersects $J\cup(J\dot+m)$, then it is contained in a dyadic subcube $K$ of $J$ or $J\dot+m$ of sidelength $\ell(K)=2^{-1}\ell(J)\geq 2^M\ell(I)\geq 2^r\ell(I)$. Since $I$ is good,
\begin{align*}
  \dist(I,K^c)>
  \ell(I)^{\gamma}\ell(K)^{1-\gamma}
  \geq 2^{M(1-\gamma)}\ell(I)\geq\abs{m}\ell(I),
\end{align*}
and hence $\dist(I\dot+m,K^c)\geq\dist(I,K^c)-\abs{m}\ell(I)>0$. This means that also $I\dot+m$ is contained in $K$, as we wanted.

We can hence decompose $\mathscr{D}_{\good}$ into collections of pairwise $m$-compatible cubes by setting
\begin{align*}
  \mathscr{D}^m_{k,v}
  :=\{I\in\mathscr{D}_{\good}:a(I)=i,b(I)=v\},\qquad
  i=0,\ldots,M(m),\quad v=0,1.
\end{align*}
The total number of these collections is $2(1+M(m))\lesssim(1+\log^+\abs{m})$.

The estimate for $A_{\good}$ now finally begins with
\begin{align*}
  \abs{A_{\good}}
  &\leq\sum_{m\in\Z^n}\sum_{I\in\mathscr{D}_{\good}}
     \sum_{\eta,\theta\in\{0,1\}^n\setminus\{0\}}
     \abs{\pair{g}{h^{\eta}_{I\dot+m}}\pair{h^{\eta}_{I\dot+m}}{Th^{\theta}_I}
     \pair{h^{\theta}_I}{f}} \\
  &\leq\Norm{g}{p'} \sum_{m\in\Z^n}\sum_{\eta,\theta}
    \sum_{k,v}\BNorm{\sum_{I\in\mathscr{D}_{k,v}^m}
      \zeta^{\theta,\eta}_{I,m}h^{\eta}_{I\dot+m}
     \pair{h^{\eta}_{I\dot+m}}{Th^{\theta}_I}
     \pair{h^{\theta}_I}{f}}{p},
\end{align*}
where the $\zeta^{\theta,\eta}_{I,m}$ are angular factors of the quantities inside the absolute values on the previous line.

For a fixed $\mathscr{D}^m_{k,v}$, the pairwise $m$-compatibility of its cubes ensures that the $d^{\eta,\theta}_{I,m,k}$ defined above, for $I\in\mathscr{D}_{k,v}^m$, form a martingale difference sequence with respect to their generated filtration, when ordered primarily according to decreasing $\ell(I)$, arbitrarily among the intervals $I$ of the same sidelength, and secondarily according to the parameter $u=0,1$. Then the support of any given $d^{\eta,\theta}_{I,m,0}$ is entirely contained in a set where all the previous members of the sequence are constant, so that its vanishing integral ensures that it is indeed a legitimate next member of a martingale difference sequence. And clearly $d^{\eta,\theta}_{I,m,1}$ has a vanishing integral separately on all the sets where $d^{\eta,\theta}_{I,m,0}$ (or the previous members of the sequence) take a given value.

Hence
\begin{align*}
  \sum_{I\in\mathscr{D}^m_{k,v}}\zeta^{\theta,\eta}_{I,m}h^{\eta}_{I\dot+m}
  \pair{h^{\eta}_{I\dot+m}}{Th^{\theta}_I}\pair{h^{\theta}_I}{f}
\end{align*}
is a martingale transforms of
\begin{align*}  
  \sum_{I\in\mathscr{D}^m_{k,v}}h^{\theta}_{I}\pair{h^{\theta}_I}{f}
  \quad\text{by the multiplying numbers}\quad
  \pm\zeta^{\theta,\eta}_{I,m}\pair{h^{\eta}_{I\dot+m}}{Th^{\theta}_I}
\end{align*}
all of which are bounded by $(1+\abs{m})^{-n-\alpha}$.

By a direct application of the defining property of UMD spaces, it then follows that
\begin{align*}
   \abs{A_{\good}}
  &\lesssim\Norm{g}{p'} \sum_{m\in\Z^n}\sum_{\eta,\theta}
    \sum_{k,v}(1+\abs{m})^{-n-\alpha}\beta_{p,X}\BNorm{\sum_{I\in\mathscr{D}^m_{k,v}}
     h^{\theta}_I\pair{h^{\theta}_I}{f}}{p},
\end{align*}
Another application of UMD with the transforming sequence of zeros and ones gives
\begin{align*}
    \BNorm{\sum_{I\in\mathscr{D}^m_{k,v}}h^{\theta}_I\pair{h^{\theta}_I}{f}}{p}
    \lesssim\beta_{p,X}\BNorm{\sum_{I\in\mathscr{D}}
      \sum_{\eta\in\{0,1\}^n\setminus\{0\}}
    h^{\eta}_I\pair{h^{\eta}_I}{f} }{p}
    =\beta_{p,X}\Norm{f}{p},
\end{align*}
and hence
\begin{align*}
   \abs{A_{\good}}
  &\lesssim\Norm{g}{p'} \sum_{m\in\Z^n}\sum_{\eta,\theta}
    \sum_{k,v}(1+\abs{m})^{-n-\alpha}\beta_{p,X}^2\Norm{f}{p} \\
  &\lesssim\Norm{g}{p'} \sum_{m\in\Z^n}(1+\log\abs{m})
    (1+\abs{m})^{-n-\alpha}\beta_{p,X}^2\Norm{f}{p}
    \lesssim\beta_{p,X}^2\Norm{g}{p'}\Norm{f}{p},
\end{align*}
and this completes the estimate for $A_{\good}$.

The considerations for $B^0_{\good}$ are almost the same, we only need to realise the $h^0_{I\dot+m}-h^0_I$ as martingale transforms of the corresponding $h^{\theta}_I$. This is achieved by setting
\begin{align*}
  d_{I,m,0}^{0,\theta}
 :=\frac{1}{3}\big(h^0_{I\dot+m}+(h^{\theta}_I)^+\big)
     -(h^{\theta}_I)^-,\qquad
  d_{I,m,1}^{0,\theta}
  :=\frac{1}{3}\big(-h^0_{I\dot+m}+2(h^{\theta}_I)^+\big),
\end{align*}
where $h^{\theta}_I=(h^{\theta}_I)^+-(h^{\theta}_I)^-$ is the splitting into positive and negative parts, so that, observing the identity $h^{0}_I=(h^{\theta}_I)^+ +(h^{\theta}_I)^-$,
\begin{align*}
 h^{\theta}_I
 =d_{I,m,0}^{0,\theta}+d_{I,m,1}^{0,\theta},\qquad
 h^0_{I\dot+m}-h^0_I
 =d_{I,m,0}^{0,\theta}-2d_{I,m,1}^{0,\theta},
\end{align*}
and hence
\begin{align*}
  \sum_{I\in\mathscr{D}_{i,j}}\zeta^{\theta,0}_{I,m}(h^{0}_{I\dot+m}-h^0_I)
  \pair{h^{0}_{I\dot+m}}{Th^{\theta}_I}\pair{h^{\theta}_I}{f}
\end{align*}
is a martingale transforms of
\begin{align*}  
  \sum_{I\in\mathscr{D}_{i,j}}h^{\theta}_{I}\pair{h^{\theta}_I}{f}
  \quad\text{by the multiplying numbers}\quad
  \{1,-2\}\cdot\zeta^{\theta,0}_{I,m}\pair{h^{0}_{I\dot+m}}{Th^{\theta}_I}
\end{align*}
all of which are bounded by $(1+\abs{m})^{-n-\alpha}$.

With obvious changes the same computation as for $A_{\good}$ then gives
\begin{align*}
   \abs{B^0_{\good}}
    \lesssim\beta_{p,X}^2\Norm{g}{p'}\Norm{f}{p}
\end{align*}
and the argument for $C^0_{\good}$, as mentioned, is dual to this. This completes the proof of Theorem~\ref{thm:Figiel}.\label{p:endFigiel}

An inspection of the argument gives the following slight generalisation:

\begin{corollary}[Figiel \cite{Figiel:90}]
Let $X$ be a UMD space and $p\in(1,\infty)$.
Let $T$ be a linear operator which satisfies  the vanishing paraproduct conditions $T(1)=T^*(1)=0$ and
\begin{align*}
  \sum_{m\in\Z^n}(1+\log^+\abs{m})\sup_I\abs{\pair{h_{I\dot+m}^{\eta}}{Th_I^{\theta}}}
  \leq C,  \qquad(\eta,\theta)\in\{0,1\}^{2n}\setminus\{(0,0)\},
\end{align*}
where the supremum is taken over all cubes $I$ in $\R^n$.
Then $T$ extends to a bounded linear operator on $L^p(\R^n;X)$, and more precisely
\begin{equation*}
  \Norm{T}{\bddlin(L^p(\R^n;X))}
  \leq C\beta_{p,X}^2,
\end{equation*}
where $C$ depends only on the dimension $n$ and the constant $C$ in the bound for the Haar coefficients.
\end{corollary}

Indeed, the bound $\abs{\pair{h_{I\dot+m}^{\eta}}{Th_I^{\theta}}}\leq C(1+\abs{m})^{-n-\alpha}$ was only used to ensure the summability of the series as in the statement of the corollary.

\section{Bourgain's inequality for translations}

This final section deals with the following result of Bourgain, which has become a cornerstone of Fourier analysis in UMD spaces. Besides Bourgain's original application \cite{Bourgain:86} to vector-valued singular integrals of convolution type, it is a key ingredient of the results in \cite{GirWei:03,HytPor:08,HytWei:06} concerning vector-valued Fourier multipliers, pseudodifferential operators, and the $T(1)$ theorem, respectively, and also in other papers. Let $\hat{f}=\mathscr{F}f$ denote the Fourier transform of $f$.

\begin{theorem}[Bourgain \cite{Bourgain:86}]\label{thm:Bourgain}
Let $X$ be a UMD space and $p\in(1,\infty)$. Let $f_j\in L^p(\R^n;X)$ be functions with $\supp\hat{f}_j\subseteq B(0,2^{-j})$. Then
\begin{equation*}
  \BNorm{\sum_j\radem_j f_j(\cdot-2^j y_j)}{p}
  \leq C\alpha_{p,X}^n\beta_{p,X}^2
     (1+\sup_j\log^+\abs{y_j})\BNorm{\sum_j\radem_j f_j}{p},
\end{equation*}
where the $\radem_j$ are independent symmetric random signs on a probability space $\Omega$, the norms are those of the space $L^p(\Omega\times\R^n;X)$, the constant $\alpha_{p,X}$ is the norm of the Hilbert transform on $L^p(\R;X)$, and $C$ is a constant depending only on the dimension $n$.
\end{theorem}

This was originally proven by Bourgain for periodic functions $f_j\in L^p(\T;X)$. It was transfered to $L^p(\R^n;X)$ by Girardi and Weis \cite{GirWei:03}, but only under the additional condition that the $y_j$ lie on the same line through the origin (which of course is no restriction when $n=1$). However, in the abovementioned applications one needs the case when $y_j\equiv y$, so the restricted statement (which can be deduced from the one-dimensional version as explained by Girardi and Weis, and only requires the constant $\alpha_{p,X}$ in place of $\alpha_{p,X}^n$) is more than sufficient.

Of course $\alpha_{p,X}\leq C\beta_{p,X}^2$ by Theorem~\ref{thm:Figiel}, and it is known by other methods (in fact, Burkholder's original proof \cite{Burkholder:83}) that one can take $C=1$ here, but as the precise connection between $\alpha_{p,X}$ and $\beta_{p,X}$ remains unknown (see \cite{GMS:10} for an up-to-date discussion), it seems better to use both constants in the estimate in the way in which they naturally appear from the proof.

We begin by deriving a dyadic analogue of this inequality, whose proof is greatly simplified by restricting ourselves to good cubes only. The estimate is the following:
\begin{equation*}
  \BNorm{\sum_j\radem_j\sum_{I\in\mathscr{D}^{\good}_j}
    h^0_{I\dot+m_j}\pair{h^0_I}{f_j}}{p}
  \lesssim\beta_{p,X}^2(1+\sup_j\log^+\abs{m_j})\BNorm{\sum_j\radem_j f_j}{p},
\end{equation*}
where $m_j\in\Z^n$ for each $j$, and $\mathscr{D}^{\good}_j:=\mathscr{D}_{\good}\cap\mathscr{D}_j$.

As before, we decompose $\mathscr{D}_{\good}$ into collections of pairwise compatible cubes. The notion of compatibility needs only slight tuning due to the fact that the relevant shifts of the cubes, $\psi:I\mapsto\psi(I):= I\dot+m_j$, $j=-\log_2\ell(I)$, now depend on the the sidelength of the cube in a more general way than before. It is now required that $I\cup\psi(I)$ and $J\cup\psi(J)$ are either disjoint or one of them, say $I\cup\psi(I)$, is contained in a dyadic subcube of either $J$ or $\psi(J)$. We let $M(\psi)$ be defined like $M(m)$ above, using $\sup_j\abs{m_j}$ in place of $\abs{m}$, and then the function $a(I)$ has exactly the same definition as before, and $b(I)$ is obviously defined by using the orbits of $\psi$ now. This provides us with a partition of $\mathscr{D}_{\good}$ into the $O(1+\sup_j\log^+\abs{m_j})$ subcollections $\mathscr{D}_{k,v}$ with $k=0,\ldots,M(\psi)$, $v=0,1$. We write $\mathscr{D}^{k,v}_j:=\mathscr{D}_{k,v}\cap\mathscr{D}_j$.

We then turn to the martingale differences, which are now functions on the product space $\Omega\times\R^n$. For each $I\in\mathscr{D}^{\good}_j$, let
\begin{equation*}
  d_{I,u}:=\radem_j(h^0_I+(-1)^uh^0_{I\dot+m_j}),\qquad u=0,1,
\end{equation*}
so that
\begin{equation*}
  \radem_j h^0_I=d_{I,0}+d_{I,1},\qquad
  \radem_j h^0_{I\dot+m_j}=d_{I,0}-d_{I,1},
\end{equation*}
and hence
\begin{equation*}
  \sum_j\radem_j\sum_{I\in\mathscr{D}^{k,v}_j}
    h^0_{I\dot+m_j}\pair{h^0_I}{f_j}
\end{equation*}
is a martingale transform of
\begin{equation*}
  \sum_j\radem_j\sum_{I\in\mathscr{D}^{k,v}_j}
    h^0_I\pair{h^0_I}{f_j}\quad
    \text{by the multiplying numbers}\quad\pm 1.
\end{equation*}
It follows that
\begin{align*}
  &\BNorm{\sum_j\radem_j\sum_{I\in\mathscr{D}^{\good}_j}
    h^0_{I\dot+m_j}\pair{h^0_I}{f_j}}{p} \\
  &\leq\sum_{k,v}\BNorm{\sum_j\radem_j\sum_{I\in\mathscr{D}^{k,v}_j}
    h^0_{I\dot+m_j}\pair{h^0_I}{f_j}}{p} \\
  &\leq\sum_{k,v}\beta_{p,X}\BNorm{\sum_j\radem_j\sum_{I\in\mathscr{D}^{k,v}_j}
    h^0_{I}\pair{h^0_I}{f_j}}{p} \\
  &\lesssim(1+\sup_j\log^+\abs{m_j})\beta_{p,X}
    \BNorm{\sum_j\radem_j\sum_{I\in\mathscr{D}_j}
     h^0_{I}\pair{h^0_I}{f_j}}{p},
\end{align*}
where the last estimate also used the contraction property of the $\radem_j$-randomised sums (pointwise in $x\in\R^n$) to return back to the full collection $\mathscr{D}_j$. Here
\begin{align*}
    \BNorm{\sum_j\radem_j\sum_{I\in\mathscr{D}_j}
     h^0_{I}\pair{h^0_I}{f_j}}{p}
     =\BNorm{\sum_j\radem_j\Exp_j f_j}{p}
     \leq\beta_{p,X}\BNorm{\sum_j\radem_j f_j}{p}.
\end{align*}
The last estimate is the vector-valued Stein inequality, which is also due to Bourgain~\cite{Bourgain:86}. It can be proven in a couple of lines directly from the definition of UMD (see~\cite{FigWoj:01}), so no difficulties of the proof are hidden into this estimate. The dyadic analogue of Bourgain's inequality has now been proven.

The next task is to compute the average
\begin{equation*}
  \Exp_{\beta}\sum_j\radem_j\sum_{I\in\mathscr{D}^{\beta}_{\good,j}}
    h^0_{I\dot+m_j}\pair{h^0_I}{f_j}
   =\pi_{\good}\cdot\Exp_{\beta}\sum_j\radem_j\sum_{I\in\mathscr{D}^{\beta}_{j}}
    h^0_{I\dot+m_j}\pair{h^0_I}{f_j}
\end{equation*}
over the random choice of $\mathscr{D}^{\beta}$. Of course, this will also satisfy the same norm bound, which all these expressions with a fixed $\beta$ have.

Recalling that $\mathscr{D}_j^{\beta}=\mathscr{D}_j^{0}+\sum_{i>j}2^{-i}\beta_i$, where the last binary expansion is uniformly distributed over $[0,2^{-j})^n$ under the random choice of $\beta$, it follows that
\begin{align*}
  &\Exp_{\beta}\sum_{I\in\mathscr{D}^{\beta}_{j}}
    h^0_{I\dot+m_j}\pair{h^0_I}{f_j} \\
  &=\int_{[0,1)^n}\sum_{k\in\Z^n}
      h^0_{2^{-j}([0,1)^n+k+u+m_j)}\pair{h^0_{2^{-j}([0,1)^n+k+u)}}{f_j}\ud u \\
  &=\int_{\R^n}
    2^{jn}h^0(2^j\cdot-u-m_j)
    \pair{h^0(2^j\cdot-u)}{f_j}\ud u\qquad\big(h^0:=h^0_{[0,1)^n}\big) \\
  &=(\varphi_{2^j}*f_j)(\cdot-2^{-j}m_j),\qquad
      \varphi(x):=\int_{\R^n}h^0(x+u)h^0(u)\ud u.
\end{align*}
Hence, as the second intermediate estimate towards Theorem~\ref{thm:Bourgain}, we obtain
\begin{align*}
  \BNorm{\sum_j\radem_j(\varphi_{2^j}*f_j)(\cdot-2^{-j}m_j)}{p}
  &\lesssim\Exp_{\beta}\BNorm{\sum_j\radem_j
  \sum_{I\in\mathscr{D}^{\beta}_{\good,j}}h^0_{I\dot+m_j}\pair{h^0_I}{f_j}}{p} \\
  &\lesssim\beta_{p,X}^2(1+\sup_j\log^+\abs{m_j})\BNorm{\sum_j\radem_j f_j}{p}.
\end{align*}
This still deviates from the final goal in two respects: there is additional smoothing on the left, and the shifts $m_j$ are restricted to integer values. Both these drawbacks may be corrected simultaneously by a Fourier multiplier technique going back to Bourgain's original proof.

By a simple change of variable, we may assume that $\supp\hat{f}_j\subseteq B(0,2^{-j-1})$ rather than just $B(0,2^{-j})$. Let $m_j$ be the integer point nearest to $y_j$, so that $z_j:=y_j-m_j\in[2^{-1},2^{-1})^n$. The quantity to be estimated has the Fourier transform
\begin{align*}
  \sum_j\radem_j \exp(-i2\pi 2^j y_j\cdot\xi)\hat{f}_j(\xi),
\end{align*}
whereas the one appearing in the intermediate estimate has transform
\begin{align*}
  \sum_j\radem_j \exp(-i2\pi 2^j m_j\cdot\xi)\hat{\varphi}(2^j\xi)\hat{f}_j(\xi),
\end{align*}
where $\hat\varphi$ is immediately computed as
\begin{align*}
  \hat\varphi(\xi)=\prod_{i=1}^n\sinc^2(\pi\xi_i),\qquad
  \sinc u:=\frac{\sin u}{u}.
\end{align*}
This function is smooth and bounded away from zero in $B(0,2^{-1})$, so that it is easy to find a function $\chi\in C_c^{\infty}(B(0,1))$ so that $\chi\hat{\varphi}\equiv 1$ in $B(0,2^{-1})$. By the support property of $\hat{f}_j$, this implies $(\chi\hat{\varphi})(2^j\xi)\hat{f}_j(\xi)\equiv\hat{f}_j(\xi)$, and hence
\begin{align*}
  &\exp(-i2\pi 2^j y_j\cdot\xi)\hat{f}_j(\xi) \\
  &=\exp(-i2\pi 2^j z_j\cdot\xi)\chi(2^j\xi)
    \times\exp(-i2\pi 2^j m_j\cdot\xi)\hat{\varphi}(2^j\xi)\hat{f}_j(\xi),
\end{align*}
The multipliers $\sigma_j(\xi):=\exp(-i2\pi z_j\cdot\xi)\chi(\xi)$, $z_j\in[-2^{-1},2^{-1})^n$, and hence their dilations $\sigma_j(2^j\xi)$ appearing above, have uniformly bounded variation in the sense that
\begin{align*}
  \int_{\R^n}\abs{\partial_1\cdots\partial_n\sigma_j(\xi)}\ud\xi\leq C,
\end{align*}
and hence the corresponding Fourier multiplier operators $T_j:=\mathscr{F}^{-1}\sigma_j(2^j\cdot)\mathscr{F}$ satisfy the estimate
\begin{align*}
  \BNorm{\sum_j\radem_j T_j g_j}{p}
  \lesssim\alpha_{p,X}^n\BNorm{\sum_j\radem_j g_j}{p}
\end{align*}
for all $g_j\in L^p(\R^n;X)$. This is another lemma of Bourgain \cite{Bourgain:86}, but not a particurly difficult one; it is based on representation of such operators as convex combinations of frequency modulations of $n$-fold products of Hilbert transforms in all coordinate directions.

The proof is hence completed by combining the estimate
\begin{align*}
  \BNorm{\sum_j\radem_j f_j(\cdot-2^{-j}y_j)}{p}
  &=\BNorm{\sum_j\radem_j T_j[\varphi_{2^j}*f_j(\cdot-2^{-j}m_j)]}{p} \\
  &\lesssim\alpha_{p,X}^n\BNorm{\sum_j\radem_j \varphi_{2^j}*f_j(\cdot-2^{-j}m_j)}{p}
\end{align*}
with the intermediate inequality already established.

\bibliography{umd-literature}
\bibliographystyle{plain}

\end{document}